\documentclass[a4paper,11pt]{article}
\usepackage{inputenc}
\usepackage{graphicx}
\usepackage{amssymb}
\usepackage{latexsym, bm}
\usepackage{multicol}
\usepackage{indentfirst}
\usepackage{amssymb,amsfonts}
\usepackage{amsmath}
\usepackage{setspace}
\usepackage{enumerate}
\usepackage{CJK}
\newtheorem{theorem}{Theorem}
\newtheorem{claim}{Claim}
\newtheorem{lemma}{Lemma}
\newtheorem{corollary}{Corollary}
\newtheorem{proposition}{Proposition}

\newtheorem{problem}{Problem}

\newtheorem{fact}{Fact}

\usepackage{amssymb}

\date{}
\begin{spacing}{1}

\begin{document}

\title{Ramsey numbers of quadrilateral versus books}

\author{Tianyu Li,\footnote{Center for Discrete Mathematics, Fuzhou University,
Fuzhou, 350108, P.~R.~China. Email: {\tt lty4765382@126.com.} } \;\;
Qizhong Lin,\footnote{Corresponding author. Center for Discrete Mathematics, Fuzhou University,
Fuzhou, 350108, P.~R.~China. Email: {\tt linqizhong@fzu.edu.cn}. Supported in part by NSFC (No.\ 1217010182).}
~~
 {Xing Peng\footnote{Center for Pure Mathematics, School of Mathematical Sciences, Anhui University, Hefei 230601, P.~R.~China. Email: {\tt x2peng@ahu.edu.cn}. Supported by the NSFC grant
(No.\ 12071002) and a Start--up Fund from Anhui University.}}
}
\maketitle

\begin{abstract}
 A book $B_n$ is a graph which consists of $n$ triangles sharing a common edge.
 In this paper, we study  Ramsey numbers of quadrilateral versus books.
 Previous results give the exact value of $r(C_4,B_n)$ for $1\le n\le 14$.
 We aim to show the  exact value of $r(C_4,B_n)$ for infinitely many $n$.
 To achieve this, we first prove that $r(C_4,B_{(m-1)^2+(t-2)})\le m^2+t$ for $m\ge4$ and $0 \leq t \leq m-1$. This improves upon a result  by Faudree, Rousseau and Sheehan (1978) which states that
\begin{align*}
r(C_4,B_n)\le g(g(n)), \;\;\text{where}\;\;g(n)=n+\lfloor\sqrt{n-1}\rfloor+2.
\end{align*}
Combining the new upper bound and constructions of $C_4$-free graphs, we are able to determine the exact value of $r(C_4,B_n)$ for infinitely many $n$.
As a special case, we show $r(C_4,B_{q^2-q-2}) = q^2+q-1$ for all prime power $q\ge4$.


\end{abstract}

\section{Introduction}
For two graphs $H_1$ and $H_2$, the {\it Ramsey number} $r(H_1,H_2)$ is the smallest integer $N$ such that for any graph $G$ on $N$ vertices, either $G$ contains $H_1$ as a subgraph or the complement $\overline{G}$ contains $H_2$ as a subgraph. Let $B_n$ be the book graph which consists of $n$ triangles sharing a common edge,
and let $C_m$ be a cycle of length $m$.

Book--cycle Ramsey numbers have attracted a great deal of attention, see  \cite{FRS78,rs,frs,shi,ll,l-p}. 
However, there is little progress towards determining the exact value of $r(C_m,B_n)$ for $m$ being a fixed even number.
As pointed out by Faudree, Rousseau and Sheehan \cite{frs}: ``we know practically nothing about $r(C_m,B_n)$ when $m$ (fixed integer) is even and greater than 4''.
Even for the case where $m=4$, we only know the exact value of $r(C_4,B_n)$ for $1\le n\le 14$.
  \begin{center}
  \begin{tabular}{c|cccccccccccccc}
  \hline
  $n$ & 1 & 2 & 3 & 4 & 5 & 6 & 7 & 8 & 9 & 10 & 11 & 12 & 13 & 14 \\

  \hline
  $r(C_4,B_n)$ & 7 & 7 & 9 & 11 & 12 & 13 & 16 & 17 & 18 & 19 & 20 & 21 & 22 & 24 \\
  \hline
  \end{tabular}

 \vspace{5mm}

   {\bf Table 1.} \; \; The value of  $r(C_4,B_n)$ for $1\leq n\leq 14$.
  \end{center}
The exact value of $r(C_4,B_1)$ was  observed by Chartrand and Schuster \cite{GS71}. For $2\leq n\leq 7$, one can see the paper by Faudree, Rousseau and Sheehan  \cite{FRS78}. For $8\leq n\leq 14$, see Tse \cite{Tse1,Tse2} in which the lower bounds were obtained by computer algorithms.


In \cite{FRS78}, the authors proved the following general upper bound for $r(C_4,B_n)$:
\begin{align*}
r(C_4,B_n)\le g(g(n)), \;\;\text{where}\;\;g(n)=n+\lfloor\sqrt{n-1}\rfloor+2.
\end{align*}
This implies that $ r(C_4,B_n)\le n+2\lfloor\sqrt{n}\rfloor+5.$ We will show the following general lower bound and the proof  will be given in the Appendix.
\begin{lemma} \label{glb}
For all large $n$, we have
$
r(C_4,B_n) \geq n+2\lfloor\sqrt{n}-6n^{0.2625}\rfloor.
$
\end{lemma}

As a special case where $q$ is a prime power, it was proved in \cite{FRS78} that
\begin{align*}
q^2+q+2\leq r(C_4,B_{q^2-q+1})\leq q^2+q+4.
\end{align*}
Let $\mathcal{G}(q)$ be the set of  graphs $G$ on $q^2+q+3$ vertices such that $G$ contains no $C_4$ and its complement contains no $B_{q^2-q+1}$. If $\mathcal{G}(q)\neq\emptyset$, then we have $r(C_4,B_{q^2-q+1})= q^2+q+4$ for prime power $q>3$. Faudree, Rousseau and Sheehan \cite{FRS78} proved that $\mathcal{G}(3) \not =\emptyset$ and proposed the following problem.
\begin{problem}[Faudree, Rousseau and Sheehan  \cite{FRS78}]
For each prime power $q>3$, whether $\mathcal{G}(q)\neq\emptyset$.
\end{problem}

The answer to Problem 1 is negative for $q=4$ as $r(C_4,B_{13})=22$, see \cite{Tse2}, which implies  $\mathcal{G}(4)=\emptyset$. The problem is still open for $q\ge5$.

\medskip
In this paper, we aim to establish the exact value of $r(C_4,B_n)$ for infinitely many $n$. At first, we will prove the following upper bound for $r(C_4,B_n)$.
\begin{theorem}\label{m-up-bou}
For any integer $m\ge4$ and $0\le t\le m-1$,
\[
r(C_4,B_{(m-1)^2+(t-2)})\le m^2+t.
\]
\end{theorem}

\medskip
\noindent
{\bf Remark.} The upper bound in Theorem \ref{m-up-bou} improves that in \cite{FRS78} by one for $m\ge4$ and $3\le t\le m-1$.
Let us point out that the upper bound gives the true value for $(m,t)=(4,0)$, and $m=3$ and $0\le t\le 2$, while it is not the case for $(m,t)=(2,2)$ as $r(C_4,B_1)=7$, see Table 1.

\medskip
For even prime power $q$, we show that the upper bound in Theorem \ref{m-up-bou} is tight by using a variation of the Erd\H{o}s-R\'enyi orthogonal polarity graph \cite{ER66,Br66}.
\begin{theorem}\label{low-bou}
Let $q\geq 4$ be an even prime power. If $0\le t \le q-1$ and $t\ne 1$, then
\[
  r(C_4,B_{(q-1)^2+(t-2)})= q^2+t.
\]
\end{theorem}


For odd prime power $q$,  we show that the upper bound in Theorem \ref{m-up-bou} is attainable by applying the $C_4$-free graphs constructed in \cite{ZC17}.
\begin{theorem}\label{odd}
Let $q\ge5$ be an odd prime power. If $q\equiv3\pmod4$, $\frac{q+1}{2}\le t\le q-1$ and $t\neq \frac{q+3}{2}$; or $q\equiv1\pmod4$, $\frac{q-1}{2}\le t\le q-1$ and $t\neq \frac{q+1}{2}$, then
\[
  r(C_4,B_{(q-1)^2+(t-2)}) = q^2+t.
\]
\end{theorem}

\medskip
As a special case when $t=q-1$, Theorem \ref{low-bou} and  Theorem \ref{odd} imply the following result.
\begin{corollary}
For all prime power $q\ge4$,
$
r(C_4,B_{q^2-q-2}) = q^2+q-1.
$
\end{corollary}

\noindent

We will use the following notations throughout the paper. Let $G=(V,E)$ be a simple graph.  We write $H\subseteq G$ if $H$ is a subgraph of $G$. 
For two disjoint subsets $X$ and $Y$, we use $E_G(X,Y)$ to denote the set of edges between $X$ and $Y$, and $e_G(X,Y)=|E_G(X,Y)|$ is the number of edges between $X$ and $Y$. Similarly, $E_G(X)$ is the set of edges contained in $X$, and $e_G(X)=|E_G(X)|$. We will write $uv$ for an edge in $G$.
For any vertex $v\in V$, let $N_G(v)$ be the set of all neighbors of $v$ and the {\it closed neighborhood} $N_G[v]$ is the set $N_G(v)\cup\{v\}$. Let $d_G(v)$ denote the degree of $v$, i.e., $d_G(v)=|N_G(v)|$. The maximum degree and the minimum degree of $G$ are denoted by $\Delta(G)$ and $\delta(G)$, respectively.
We always omit the subscript $G$ whenever it is clear under the context.

The rest of the paper is organized as follows. In Section 2, we will present the proof of Theorem \ref{m-up-bou}. In Section 3, we will prove Theorem \ref{low-bou} and Theorem  \ref{odd}. Finally, we will mention a few concluding remarks in Section 4.


\section{Proof of Theorem \ref{m-up-bou}}
A classical result due to Parsons \cite{TP75} gives exact value of $r(C_4,K_{1,n})$ for infinitely many $n$. The author also proved the following upper bound.

\begin{lemma}[Parsons \cite{TP75}]\label{pars}
For any integer $n\geq 2$,
 $r(C_4,K_{1,n})\leq n+\lfloor\sqrt{n-1}\rfloor+2$.
 Moreover, if $n=k^2+1$ and $k\ge1$, then
 $r(C_4,K_{1,n})\leq n+\lfloor\sqrt{n-1}\rfloor+1$.
\end{lemma}


The following proposition can be easily checked by using the double-counting method.
\begin{proposition}\label{C4-eg}
Let $G$ be a graph. If $\sum_{v\in V(G)} {d(v)\choose 2}>{|V(G)| \choose 2}$,
then $C_4\subseteq G$.
\end{proposition}


\noindent
{\bf Assumption:} \ Let $N=m^2+t$, where $m\ge4$ and $0 \leq t \leq m-1$.
Suppose $G$ is a graph on $N$ vertices such that
$$C_4 \nsubseteq G,\;\; \text{and} \;\; B_{(m-1)^2+t-2} \nsubseteq \overline G.$$

We will prove some properties satisfied by the graph $G$ under the above assumption.

\begin{claim}\label{mindegree}
$\delta(G)=m$.
\end{claim}
\noindent\textbf{Proof.} If $\delta(G) \leq m-1$, then $\Delta(\overline{G})\geq N-1-\delta(G) \geq m^2-m+t$. Let $v$ be a vertex such that $d_{\overline{G}}(v) \geq m^2-m+t$. As $m^2-m+t \geq (m-1)^2+t-2+\lfloor \sqrt{(m-1)^2+t-3}\rfloor+2$, by Lemma \ref{pars}, either $N_{\overline{G}}(v)$ contains a $C_4$ or its complement contains a star $K_{1,(m-1)^2+t-2}$. In the former case, there exists a $C_4$ in $G$, a contradiction. In the latter case, the star together with the vertex $v$ form a book $B_{(m-1)^2+t-2}$ in $\overline{G}$, a contradiction. It follows that $\delta(G) \geq m$. If $\delta(G) \geq m+1$, then it is clear that $\sum_{v\in V(G)} {d(v)\choose 2} \geq N \binom{m+1}{2}> \binom{N}{2}$.
Therefore, we get that $C_4 \subseteq G$ by Proposition \ref{C4-eg}, a contradiction and the claim follows. \hfill    $\Box$

\begin{claim}\label{m-adj}
 Let $u$ and $w$ be two vertices of degree $m$ in $G$. If $|N(u)\cap N(w)|=1$, then $uw\in E(G)$.
\end{claim}

\noindent\textbf{Proof.} Suppose that $u$ and $w$ are not adjacent in $G$. Since $|N(u)\cap N(w)|=1$, it follows that $u$ and $w$ have $m^2+t-2-(2m-1)=(m-1)^2+(t-2)$ common non-neighbors, which implies that $B_{(m-1)^2+(t-2)}\subseteq \overline{G}$, a contradiction.
\hfill$\Box$

\medskip
Let  $v \in V(G)$ be a vertex of degree $m$ and  $N(v)=\{v_1,v_2,\dots,v_m\}$. For each $1\le i\le m$, we define (see Fig. 1)
\begin{align*}
&A_i=N(v_i)\setminus N[v], \; \; \text{where $N[v]$ is the closed neighborhood of $v$ in $G$},\\
&B_i=(\cup_{x\in A_i}N(x))\setminus [\{v_i\}\cup (\cup_{1\le i\le m}A_i)],\\
&B=V(G) \setminus [N[v] \cup (\cup_{i=1}^m A_i)].
\end{align*}

\begin{figure}
\begin{center} \label{aibi}
\includegraphics[scale=0.25]{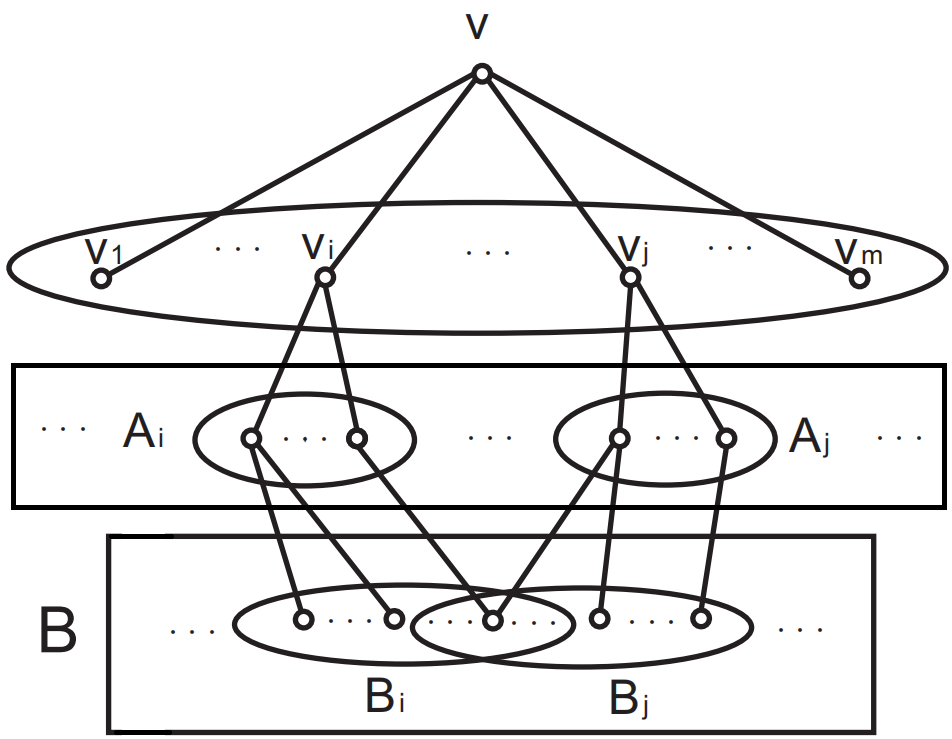}
\end{center}
\begin{center}
{\bf Fig. 1}: \; The structure of $G$.
\end{center}
\end{figure}

The following claim is clear from the assumption $C_4\nsubseteq G$.
\begin{claim}\label{structure1}
The graph $G$ satisfies the following:

1. For any vertex $u\in V(G)$, edges in $N(u)$ form a matching.

2. For $1\le i\le m$, edges in $A_i$ form a matching.

3. For $1\le i<j\le m$, $A_i\cap A_j=\emptyset$.

4. For $1 \leq i < j \leq m$, $E(A_i, A_j)$ is a matching.

5. If $v_i$ and $v_j$ are adjacent in $G$, then $E(A_i, A_j)=\emptyset$.

6. For $1\le i\le m$, any vertex in $B_i$ has exactly one neighbor in $A_i$.
\end{claim}

\begin{claim}\label{structure2}
The graph $G$ also satisfies the following:

1. There are at most two vertices of degree $m$ in $N(v)$.

2. For $1\le i\le m$ and any vertex $x\in A_i$, it holds $d(x)\ge m+1$.
\end{claim}
\textbf{Proof.} By Claim \ref{m-adj}, any pair of vertices of degree $m$ in $N(v)$ must be adjacent. If there are three vertices of degree $m$ in $N(v)$, then we can get a $C_4$, which is a contradiction. For the second part, if there exists a vertex $x\in A_i$ with $d(x)= m$, then Claim \ref{m-adj} implies that $v$ and $x$ are adjacent by noting $v_i\in N(v)\cap N(x)$. A contradiction.
\hfill$\Box$

\medskip
We also have the following facts.

\begin{fact}\label{d-A} We have
$
\sum\limits_{i=1}^{m}d(v_i)=m+\sum\limits_{i=1}^{m}|A_i|+2e(N(v)).
$
\end{fact}

\begin{fact}\label{u-1} For $1\le i\le m$,
$
\sum\limits_{x\in A_i}(d(x)-1)=2e(A_i)+\sum\limits_{1\le j\le m,j\ne i}e(A_i,A_j)+|B_i|.
$
\end{fact}
{\bf Proof.} Observe that
$
\sum\limits_{x\in A_i}(d(x)-1)=2e(A_i)+\sum\limits_{1\le j\le m,j\ne i}e(A_i,A_j)+e(A_i,B_i).
$
Note that $e(A_i,B_i)=|B_i|$ by Part 6 of Claim \ref{structure1}. The assertion follows. \hfill$\Box$

\medskip
We can estimate the size of $B_i$ as follows.
\begin{claim}\label{Bi-Ai}
For each $v_i\in N(v)$ where $1\le i\le m$, if there exists some vertex $v_j$ (at most one such vertex) which is adjacent to $v_i$, then $|B_i| \ge |A_i|$.
\end{claim}
{\bf Proof.}  The assumption $C_4\nsubseteq G$ yields that for any vertex $x\in A_i$,

\medskip
(i) $x$ has at most one neighbor in $A_i$,

(ii) $x$ has no neighbor in $A_j$ as $v_i$ is adjacent to $v_j$,

(iii) $x$ has at most one neighbor in $A_k$ for each $1 \leq k \leq m$ and $k \neq i,j$.
\medskip

\noindent
For each vertex $x\in A_i$, note that  $d(x)\ge m+1$ by Part 2 of Claim \ref{structure2}. It follows that $x$ has at least one neighbor in $B_i$.
Therefore, we obtain $|B_i|\ge |A_i|$ as claimed since any vertex in $B_i$ has exactly one neighbor in $A_i$ by Part 6 of Claim \ref{structure1}. \hfill$\Box$

\begin{claim}\label{edges}
We have $e(N(v))\geq \lceil \frac{m-2}{2} \rceil$.
\end{claim}
{\bf Proof.} Let $\tau$ be the number of vertices of degree $m$ in $N(v)$, and $\sigma$ be the number of vertices  in $N(v)$ of degree at least $m+2$. Each of the remaining vertices in $N(v)$ has degree $m+1$. We first show $E(N(v)) \neq \emptyset$. If it is not the case, then we get $\tau <2$ as two vertices of degree $m$ with one common neighbor $v$ are adjacent by Claim \ref{m-adj}.
 Moreover, we have $|A_i|=d(v_i)-1$ for $1 \leq i \leq m$. Note that $\delta(G)=m$ from Claim \ref{mindegree} and  $A_i\cap A_j=\emptyset$ for $1\le i<j\le m$ from Part 3 of Claim \ref{structure1}. We get
\begin{align*}\label{Nv-Ai}
  |N[v] \cup A_1 \cup \cdots \cup A_m|&=m+1+\sum_{i=1}^m|A_i| \\
&= m+1+\tau(m-1)+(m-\tau-\sigma)m+\sigma(m +1) \\
 &= m^2+m+1+\sigma-\tau.
\end{align*}
As $ |N[v] \cup A_1 \cup \cdots \cup A_m| \leq N=m^2+t \leq m^2+m-1$, it follows that $\tau \geq 2$, a contradiction. Therefore, $E(N(v)) \neq \emptyset$.

Now we assume that $v_1v_2$ is an edge in $N(v)$. Since $\delta(G)=m$, it follows from Claim \ref{Bi-Ai} that
$|B| \geq |B_1| \geq |A_1|=d(v_1)-2 \geq m-2.$
By Fact \ref{d-A},
\[
\sum_{i=1}^m |A_i|=\sum_{i=1}^m d(v_i)-m-2e(N(v))\ge2m+(m-2)(m+1)-m-2e(N(v)) \geq m^2-2-2e(N(v)).
\]
Therefore,
\begin{align*}
m^2+m-1 \geq |V(G)| &=|N[v] \cup A_1 \cup \cdots \cup A_m \cup B|\\
                  & \geq m+1+(m^2-2-2e(N(v)))+m-2\\
                   &= m^2+2m-3-2e(N(v)).
\end{align*}
Consequently, $e(N(v)) \geq \tfrac{m-2}{2}$. Since $e(N(v))$ is an integer, we get $e(N(v)) \geq \lceil \frac{m-2}{2} \rceil$ as claimed no matter the parity of $m$. \hfill $\Box$

\medskip
We are now ready to prove Theorem \ref{m-up-bou}.

\medskip
\noindent\textbf{Proof of Theorem \ref{m-up-bou}.}
Let $N=m^2+t$. We will show that $r(C_4,B_{(m-1)^2+(t-2)})\le N$ for $m\ge4$ and $0\le t\le m-1$.
Suppose that  $G$ is a graph on $N$ vertices such that
 $$C_4 \nsubseteq G,\;\; \text{and} \;\; B_{(m-1)^2+t-2} \nsubseteq \overline G.$$
We will find a contradiction, which will complete the proof of Theorem \ref{m-up-bou}.

For $0\leq t\leq 2$, Lemma \ref{pars} implies that
\begin{align*}
\Delta(\overline{G})&\leq r(C_4,K_{1,(m-1)^2+(t-2)})-1
\\&\leq (m-1)^2+(t-2)+\left\lfloor\sqrt{(m-1)^2+(t-2)-1}\right\rfloor+1
\\&\leq m^2-m+t-2.
\end{align*}
For $t=3$, the second assertion of Lemma \ref{pars} implies that
\begin{align*}
\Delta(\overline{G})&\leq r(C_4,K_{1,(m-1)^2+1})-1
\leq (m-1)^2+1+\left\lfloor\sqrt{(m-1)^2}\right\rfloor
= m^2-m+1.
\end{align*}
Therefore, we conclude with $\Delta(\overline{G})\le m^2-m+t-2$ for $0\le t\le3$.
It follows that
 \[
 \delta(G)\geq N-1-\Delta(\overline{G})\ge m+1.
 \]
This is a  contradiction to Claim \ref{mindegree}.
 In the following, we assume that $4\le t\le m-1$.

Let $v$ be a vertex of degree $m$ in $G$ and $N(v)=\{v_1,v_2,\dots,v_m\}$. For $1\le i\le m$, we define sets $A_i$, $B_i$, and $B$ as before, see Fig. 1. Clearly, $B_i \subseteq B$. We will obtain a contradiction by showing $|B_i|>|B|$ for some $1 \leq i \leq m$. Let $\tau$ and $\sigma$ be defined as in Claim \ref{edges}. 
By Fact \ref{d-A},
\begin{align*}
\sum\limits_{i=1}^{m}|A_i|&=\sum\limits_{i=1}^{m}d(v_i)-m-2e(N(v))
\\&\ge\tau m+\sigma(m+2)+(m-\tau-\sigma)(m+1)-m-2e(N(v))
\\&\ge m^2+\sigma-\tau-2e(N(v)).
\end{align*}
From the definition of $A_i$, $B_i$, and $B$, we obtain that
\begin{align}\label{|B|}
\nonumber|B|=N-\left(|N[v]|+\sum\limits_{i=1}^{m}|A_i|\right)
&\nonumber\le(m^2+t)-(m+1+m^2+\sigma-\tau-2e(N(v)))
\\&= 2e(N(v))+(t-m)-1+\tau-\sigma.
\end{align}
Note that Claim \ref{edges} implies that $$e(N(v)) \geq \lceil\tfrac{m-2}{2}\rceil\ge1,$$
so we may assume that $v_1v_2$ is an edge in $N(v)$.

If $\tau=0$, i.e., $N(v)$ contains no vertex of degree $m$, then each vertex $v_i\in N(v)$ has degree at least $m+1$. For this case, $|B|\le t-1$ from (\ref{|B|}). Since $v_1v_2$ is an edge, it follows from Claim \ref{Bi-Ai} that $|B_1|\ge |A_1|=d(v_1)-2\ge m-1>|B|$, a contradiction.

 If $\tau=1$, then one of $v_1$ and $v_2$ has degree at least $m+1$.  We  assume that $d(v_2)\ge m+1$ and $d(v_\ell)=m$ for some $\ell \neq 2$. As $v_1v_2$ is an edge, Claim \ref{Bi-Ai} implies that $|B_2|\ge |A_2|\ge m-1$. Recall inequality \eqref{|B|}. We get  that $|B|<m-1 \leq |B_2|$ unless $m$ is even, $e(N(v))=\frac{m}{2}$, and $\sigma=0$.
 In this case, we have $d(v_i) = m+1$ for $1 \leq i \neq \ell \leq m$.
As we assumed that $m \ge 4$ is even and $4 \leq t \leq m-1$, we get $m \geq 6$ and $N(v)$ contains at least three edges. Thus we can assume
that $v_3v_4$ is an edge and   $d(v_3)=d(v_4)=m+1$.
 By Part 2 of Claim \ref{structure2}, $d(x)\ge m+1$ for any vertex $x\in A_3$. Moreover, $e(A_3)\le m/2-1$ since $m-1$ is odd and edges in $A_3$ induces a matching by Part 2 of Claim \ref{structure1}. Note that Part 5 of Claim \ref{structure1} implies that $e(A_3,A_4)=0$ as $v_3v_4$ is an edge.
  Part 4 of Claim \ref{structure1} gives that $e(A_\ell,A_3)\le m-2$ as $|A_\ell|=m-2$.
  Similarly, we get $e(A_3,A_k)\le m-1$ for $1\le k\le m$ and $k\not \in \{3,\ell\}$.
 Therefore, combining these facts with Fact \ref{u-1}, we obtain that
\begin{align*}
m(m-1)&\le\sum\limits_{x\in A_3}(d(x)-1)=2e(A_3)+\sum\limits_{1\le k\le m, k\ne 3}e(A_3,A_k)+|B_3|
\\&\le (m-2)+((m-2)+(m-1)(m-3))+|B_3|.
\end{align*}
Thus $|B_3|\ge m+1>t\ge |B|$, a contradiction.

In the following, we assume that $\tau=2$. Recall the assumption that $v_1v_2$ is an edge. Since two vertices of degree $m$ must be adjacent, we are good to assume that $v_1$ and $v_2$ have degree $m$.
Note that the inequality (\ref{|B|}) implies that
\begin{align}\label{|B|3}
|B|\le 2e(N(v))+(t-m)+1-\sigma\le 2e(N(v))+(t-m)+1.
\end{align}

{\bf Case 1} \; $m$ is even and $4\le t\le m-1$.

\medskip
If $e(N(v))=m/2-1$, then $|B| \leq t-1 \leq m-2$ by the inequality (\ref{|B|3}). For this case, $m \geq 6$ is even and $N(v)$ contains at least two edges. We may further assume that $v_3v_4$ is an edge in $N(v)$.
It follows by Claim \ref{Bi-Ai} that $|B_3|\ge |A_3|\ge m-1>|B|$, a contradiction.

If $e(N(v))=m/2$, then $|B|\le t+1$ by the inequality  (\ref{|B|3}).
If there exists some $3\le i\le m$ with $d(v_i)=m+1$, says $v_3$, then $|A_3|=m-1$.
A similar argument as above, we obtain that

\medskip
(i) $d(x)\ge m+1$ for any vertex $x\in A_3$,

(ii) $e(A_3)\le m/2-1$,

(iii) $e(A_3,A_4)=0$,

(iv) $e(A_k,A_3)\le m-2$ for $k=1,2$ and $e(A_3,A_k)\le m-1$ for $5\le k\le m$.

\medskip\noindent
It follows from Fact \ref{u-1} that
\begin{align*}
m(m-1)&\le\sum\limits_{x\in A_3}(d(x)-1)=2e(A_3)+\sum\limits_{1\le k\le m, k\ne 3}e(A_3,A_k)+|B_3|
\\&\le (m-2)+(2(m-2)+(m-1)(m-4))+|B_3|.
\end{align*}
Thus $|B_3| \geq m+2 >t+1 \geq |B|$, a contradiction.

Now we assume that  $d(v_i)\ge m+2$ for each $3\le i\le m$. Note that $|A_i|\ge m$ for each $3\le i\le m$. Therefore, $\sum\limits_{i=1}^{m}|A_i|\ge 2(m-2)+(m-2)m =m^2-4$. Consequently, we  can bound the size of $B$ from above as follows:
\begin{align*}
|B|=N-\left(|N[v]|+\sum\limits_{i=1}^{m}|A_i|\right)
\le m^2+t-(m+1+m^2-4)
=3+t-m.
\end{align*}
It follows from Claim \ref{Bi-Ai} that for each $3\le i\le m$,
\[
|B_i|\ge |A_i|\ge m> 3+t-m\ge |B|,
\]
which leads to a contradiction.

\medskip
{\bf Case 2} \; $m$ is odd and $4\le t\le m-1$.

\medskip
We only need  to check the case where $e(N(v))=(m-1)/2$ because of $e(N(v)) \geq \lceil\tfrac{m-2}{2}\rceil$. Then $N(v)$ contains at least two edges. We may assume that $v_1v_2$ and $v_3v_4$ are two edges.
Note that $|B|\le t$ by the inequality (\ref{|B|3}).


If $d(v_3)\ge m+2$, then $|B_3|\ge |A_3|\ge m>t\ge |B|,$
a contradiction. Thus we assume $d(v_3)= m+1$, which implies that $|A_3|=m-1$. Clearly, $e(A_3)\le (m-1)/2$ since edges in $A_3$ is a matching. A similar argument as above yields that  $d(x)\ge m+1$ for each vertex $x\in A_3$. Moreover,
$e(A_3,A_4)=0$ as $v_3v_4$ is an edge. Similarly, $e(A_k,A_3)\le m-2$ for $k \in \{1,2\}$ as  $|A_1|=|A_2|=m-2$, and $e(A_3,A_k)\le m-1$ for $5\le k\le m$. It follows from Fact \ref{u-1} that
\begin{align*}
m(m-1)&\le\sum\limits_{x\in A_3}(d(x)-1)=2e(A_3)+\sum\limits_{1\le k\le m, k\ne 3}e(A_3,A_k)+|B_3|
\\&\le (m-1)+(2(m-2)+(m-1)(m-4))+|B_3|.
\end{align*}
Thus $|B_3|\ge m+1>t\ge |B|$, a contradiction.

Therefore, there does not exist a graph $G$ on $m^2+t$ vertices such that $C_4\nsubseteq G$ and $B_{(m-1)^2+(t-2)}\nsubseteq \overline{G}$ provided $m\ge4$ and $0\le t\le m-1$.
The proof of Theorem \ref{m-up-bou} is complete now.
\hfill$\Box$

\section{Proofs of Theorem \ref{low-bou} and Theorem \ref{odd}}\label{pf-th2-3}

We recall the definition of the Erd\H{o}s-R\'{e}nyi orthogonal polarity graph.  If $q$ is a prime power, then let $F_q$ be the Galois field with $q$ elements and $F_q^*=F_q\setminus\{0\}$.
 For two vectors $(a_1,a_2,a_3)$, $(b_1,b_2,b_3)\in (F_q^3)^*$, we write $(a_1,a_2,a_3)\equiv(b_1,b_2,b_3)$ if there exists an element $\lambda\in F_q^*$ such that $(a_1,a_2,a_3)=\lambda(b_1,b_2,b_3)$.  Note that `$\equiv$' is an equivalence relation over $(F_q^3)^*=F_q^3\setminus\{(0,0,0)\}$.
 Let $\langle a_1,a_2,a_3 \rangle$ denote the equivalence class containing $(a_1,a_2,a_3)$, and let $V_q$ be the set of all equivalence classes. Obviously, $|V_q|=q^2+q+1$.

Let $G_q$ be a graph on vertex set $V_q$ in which any pair of vertices $\langle a_1,a_2,a_3 \rangle$ and $\langle b_1,b_2,b_3 \rangle$ are adjacent in $G_q$ if and only if $a_1b_1+a_2b_2+a_3b_3=0.$
Therefore, the graph $G_q$ contains a loop at the vertex $\langle a_1,a_2,a_3 \rangle$ if and only if $a_1^2+a_2^2+a_3^2=0.$
The Erd\H{o}s-R\'{e}nyi graph $ER_q$ is obtained from $G_q$ by deleting all loops. This construction can be seen in Erd\H{o}s, R\'{e}nyi and S\'{o}s \cite{ER66} (or Brown \cite{Br66}).

\medskip
We list some useful properties of $ER_q$ as follows.
\begin{lemma}[Erd\H{o}s, R\'{e}nyi and S\'{o}s \cite{ER66}]\label{er-pro}
For any prime power $q$, $ER_q$ has the following properties.

\medskip
(i) The diameter of $ER_q$ is 2;

(ii) $C_4\nsubseteq ER_q$;

(iii) For any vertex $v\in V_q$, $d(v)=q$ or $d(v)=q+1$.
\end{lemma}

Note that $\langle a_1,a_2,a_3 \rangle$ is a $q$-vertex in $ER_q$ if and only if $a_1^2+a_2^2+a_3^2=0$.
The following property is an immediate corollary.
\begin{lemma}\label{W-ind-set}
The set consisting of all vertices of degree $q$ in $ER_q$ forms an independent set.
\end{lemma}

\subsection{Proof of Theorem \ref{low-bou}}

 Let $q\ge4$ be an even prime power. Note that $\langle 1,1,1 \rangle$ is a vertex of degree $q+1$ in $ER_q$.
Let
$N(\langle 1,1,1 \rangle)=\{w_1,w_2,\dots,w_{q+1}\}$,  and $A_{w_i}=N(w_i)\setminus N[\langle 1,1,1 \rangle]$ for $1\le i\le q+1$.
Since $N(\langle 1,1,1 \rangle)$ is an independent set by Lemma \ref{W-ind-set}, we have that $|A_{w_i}|=q-1$ for $1\le i\le q+1$.

\begin{claim}\label{W-N}
We have $N(\langle 1,1,1 \rangle)$ consists of all vertices of degree $q$ in $ER_q$.
\end{claim}
{\bf Proof.} Let  $\langle a_1,a_2,a_3 \rangle$ be a vertex of degree $q$. As the characteristic of $F_q$ is 2, we have
$
0=a_1^2+a_2^2+a_3^2=(a_1+a_2+a_3)^2.
$
 It follows that $a_1+a_2+a_3=0$, which implies that $\langle a_1,a_2,a_3 \rangle$ is a neighbor of $\langle 1,1,1 \rangle$. \hfill$\Box$

\medskip
 Note that $A_{w_i}\cap A_{w_j}=\emptyset$ for $1\le i<j\le q+1$, and
$|N[\langle 1,1,1 \rangle]|+\sum_{1\le i\le q+1}|A_{w_i}|=q^2+q+1.$
 It follows that
$\{\langle 1,1,1 \rangle\}$, $N(\langle 1,1,1 \rangle)$ and $\cup_{i=1}^{q+1} A_{w_i}$ form a partition of $V_q$, i.e.,
\begin{align*}\label{partition}
V_q=\{\langle 1,1,1 \rangle\}\cup N(\langle 1,1,1 \rangle)\cup (\cup_{1\le i\le q+1}A_{w_i}).
\end{align*}
\begin{claim}\label{match}
 For $1\le i<j\le q+1$, $E(A_{w_i},A_{w_j})$ is a perfect matching in $ER_q$.
\end{claim}
{\bf Proof.} Clearly, $E(A_{w_i},A_{w_j})$ is a matching since $C_4\nsubseteq ER_q$. Note that  $N(\langle 1,1,1 \rangle)$ forms an independent set. Moreover, each vertex $x\in A_{w_i}$ must be adjacent to some vertex in $N(w_j)$ as the diameter of $ER_q$ is 2 by Lemma \ref{er-pro} (i). It follows that $x$ must be adjacent to some vertex in $A_{w_j}$ as $x$ is not adjacent to $\langle 1,1,1\rangle$ or $w_j$. The claim follows. \hfill$\Box$

\medskip

 Claim \ref{match} gives the following proposition.
\begin{claim}\label{Ai}
For $1\le i\le q+1$, $A_{w_i}$ induces an independent set in $ER_q$.
\end{claim}

We will construct a graph $H_q^t$ on $q^2+t-1$ vertices obtained from $ER_q$ such that $C_4\nsubseteq G$ and $B_{(q-1)^2+(t-2)}\nsubseteq \overline{H_q^t}$ for $0\le t \le q-1$ and $t\ne 1$.
The construction other than the case $t=0$ can be found in Parsons \cite{TP76}.
For $t=0$, let $H_q^0=ER_q\setminus N[\langle 1,1,1\rangle]$.
For $2\le t\le q-1$, let $A_{w_1}=\{u_1,u_2,\dots,u_{q-1}\}$, and we define
\[
H_q^t=ER_q\setminus\{w_1,u_{t-1},\dots,u_{q-1}\}.
\]
It is clear that $H_q^t$ has $q^2+t-1$ vertices, and $C_4\nsubseteq H_q^t$. Moreover,
\[
d_{H_q^t}(u)=\left\{
        \begin{array}{ll}
         q,\;\; &\text{$u\in\{\langle 1,1,1\rangle, w_2,\dots,w_{q+1},u_1,\dots,u_{t-2}\}$},
           \\ q, \;\; &\text{$u\in A_{w_2}\cup\cdots\cup A_{q+1}$ is adjacent to some $u_i$ for $t-1\le i\le q-1$},
          \\  q+1,\;\; &\text{otherwise}.

        \end{array}
    \right.
\]
Therefore, $H_q^t$ has minimum degree at least $q$.
\medskip

Note that $|N(u)\cap N(v)|\leq 1$ for any two vertices $u,v\in V(H_q^t)$ since $C_4\nsubseteq H_q^t$. It follows that $|N(u)\cup N(v)|\geq 2q-1$ as $\delta(H_q^t)\ge q$. Thus $u$ and $v$ have at most
\begin{align*}
|V(H_q^t)|-2-|N(u)\cup N(v)|\le (q^2+t-1)-2 -(2q-1)
=(q-1)^2+t-3
\end{align*}
common non-neighbors, i.e., $B_{(q-1)^2+(t-2)}\nsubseteq \overline{H_q^t}$.
Consequently, $r(C_4,B_{(q-1)^2+(t-2)})\ge q^2+t$ for even prime power $q\ge4$, $0\le t \le q-1$ and $t\ne 1$, which together with Theorem \ref{m-up-bou} complete the proof of Theorem \ref{low-bou}.
\hfill$\Box$

\subsection{Proof of Theorem \ref{odd}}

Let $q\ge 5$ be an odd prime power. 
Note that $\langle 0,0,1 \rangle$ is a vertex of degree $q+1$ in $ER_q$.
In the following, let $N(\langle 0,0,1 \rangle)=\{w_1,w_2,\dots,w_{q+1}\}$ and $A_{w_i}=N(w_i)\setminus N[\langle 0,0,1 \rangle]$ for $1\le i\le q+1$.

\begin{lemma}[Zhang, Chen and Cheng \cite{ZC17}]\label{per-mat}
Let $q\ge 5$ be an odd prime power. If $w_i$ is a vertex of degree $q+1$, and $w_j$ is non-adjacent to $w_i$, then $E(A_{w_i},A_{w_j})$ induces a perfect matching, where $1\le i,j\le q+1$.
\end{lemma}

The authors \cite{ZC17} obtained the following property which describes the local structures for the corresponding Erd\H{o}s-R\'enyi orthogonal polarity graph.

\medskip
(i) If $q\equiv 3\pmod4$, then $w_i$ has degree $q+1$ for $1\le i\le q+1$.
Moreover, $N(\langle 0,0,1 \rangle)$ induces a perfect matching.  We assume that $w_iw_{i+1}$ is an edge for each odd $i$ with $1\le i\le q$, and $w_i$ is adjacent to exactly two vertices of degree $q$ for $1\le i\le (q+1)/2$. Additionally,
$\{\langle 0,0,1 \rangle\}$, $N(\langle 0,0,1 \rangle)$ and $\cup_{i=1}^{q+1}A_{w_i}$  form a partition of $V_q$.

\medskip
(ii) If $q\equiv 1\pmod4$, then there are exactly two vertices of degree $q$ in $N(\langle 0,0,1 \rangle)$, say $w_q$ and $w_{q+1}$.
   In addition, $w_qw_{q+1}\not\in E(ER_q)$ and $N(\langle 0,0,1 \rangle)\setminus\{w_q, w_{q+1}\}$ induces a perfect matching. We assume that $w_iw_{i+1}$ is an edge for each odd $i$ with $1\le i\le q-2$. Furthermore, for each $1\le i\le (q-1)/2$, the vertex $w_i$ is adjacent to exactly two vertices of degree $q$. Note that  $\{\langle 0,0,1 \rangle\}$, $N(\langle 0,0,1 \rangle)$ and $\cup_{i=1}^{q+1}A_{w_i}$  form a partition of $V_q$.

\medskip
We now construct a $C_4$-free graph $G_q^t$ which is a simple variation of the one constructed in \cite{ZC17}, where $q\equiv3\pmod4$, $\frac{q+1}{2}\le t\le q-1$ and $t\neq \frac{q+3}{2}$; or $q\equiv1\pmod4$, $\frac{q-1}{2}\le t\le q-1$ and $t\neq \frac{q+1}{2}$.

Let $G^*=ER_q\setminus \{ \langle 0,0,1 \rangle, w_{t+1},\dots,w_{q+1} \}$, and we define
\[
G_q^t=\left\{
        \begin{array}{ll}
           G^*, \;\; &\text{$t$ is even},
          \\  G^*+w_{t-1}{w_t}-E(A_{w_{t-1}},A_{w_t}),\;\; &\text{$t$ is odd},

        \end{array}
    \right.
\]
i.e., for odd $t$, $G_q^t$ is obtained from $G^*$ by adding $w_{t-1}{w_t}$ and deleting edges in $E(A_{w_{t-1}},A_{w_t})$.

It is clear that $G_t$ has $q^2+t-1$ vertices, $C_4\nsubseteq G_q^t$ and $\delta (G_q^t)\ge q$.
Therefore, a similar argument as that in Theorem \ref{low-bou} gives $r(C_4,B_{(q-1)^2+(t-2)})\ge q^2+t$, which together with Theorem \ref{m-up-bou} complete the proof of Theorem \ref{odd}.\hfill$\Box$

\section{Concluding remarks}
For $(m,t)=(5,2)$ and $(m,t)=(5,4)$, Theorem \ref{odd} gives us the following corollary.
\begin{corollary}
We have $r(C_4,B_{16})=27$ and  $r(C_4,B_{18})=29$.
\end{corollary}

Recall $\mathcal{G}(q)$ is the set of  graphs $G$ on $q^2+q+3$ vertices such that $G$ contains no $C_4$ and its complement contains no $B_{q^2-q+1}$. The following problem is still open.
\begin{problem} [Faudree, Rousseau and Sheehan \cite{FRS78}]
For each prime power $q\ge5$, whether $\mathcal{G}(q)\neq\emptyset$.
\end{problem}

We mention that one cannot  improve the upper bound $r(C_4,B_{q^2-q+1})\le q^2+q+4$ (see \cite{FRS78}) by repeating the argument in
the proof of Theorem \ref{m-up-bou}.

\noindent
{\bf Appendix}

\medskip
\noindent
{\bf Proof of Lemma \ref{glb}.}  The main idea of the proof  comes from \cite{befrs}. Suppose that there exists a graph $G$ on $n + 2m$ vertices satisfying $C_4\nsubseteq G$ and $\delta(G)\ge m$. Then the lower bound $r(C_4,B_n)>n+2m$ follows by noting that any two non-adjacent vertices in such a graph have at most $(n+2m)-2-1-2(m-1)=n-1$ common non-neighbors. In the following, we will prove the existence of such graph.

 Let $\alpha=0.525$ and  $p$ be the smallest prime with $p\ge n^{1/2}+\frac12$.
  A result in \cite{bhp}  asserts that for sufficiently large $k$ there is a prime in the interval $(k,k+k^{0.525}]$.
 Thus  we get $p\le n^{1/2}+n^{\alpha/2}+1$. For each prime $p$, it is well-known \cite{ER66,Br66} that there exists a graph $ER_p$ of order $N = p^2 +p+ 1$ which contains no $C_4$ and in which the degree of each vertex is $p$ or $p + 1$.
Set $m = \lfloor n^{1/2}-6n^{\alpha/2}\rfloor$ and $d = N-(n + 2m)$.
Let us randomly delete $d$ vertices from $ER_p$ to obtain a random graph $G$. We say a vertex $v$ of degree $p$ in $ER_p$ ``bad'' if it is not deleted and has degree less than $m$ in $G$.
Let $B_v$ denote the event that $v$ is bad. Note that a vertex $v$ of degree $p$ is bad if and only if $v$ is not deleted and there are $t$  neighborhood of $v$ are removed for  some $t> p-m$. Therefore,
\[
\Pr(B_v)=\sum_{t > p-m}\frac{{p\choose t}{{N-p-1}\choose {d-t}}}{{N\choose d}}
<\sum_{t > p-m}\left(\frac{epd}{t(N-d)}\right)^t<N\left(\frac{2e}{6}+o(1)\right)^{6n^{\alpha/2}}
\]
by noting $p< n^{1/2}+n^{\alpha/2}+1$, $d= N-(n + 2m)\le (2+o(1))n^{(1+\alpha)/2}$, $t> p-m > 6n^{\alpha/2}$ and $N-d=n+2m>n$. It follows that $N\Pr(B_v)\to 0$ as $n\to\infty$. Clearly, a similar argument works for a vertex $v$ with degree $p+1$. Thus with high probability that all vertices
are good if $n$ is large. Therefore, there exists a graph $G$ on $n + 2m$ vertices satisfying $C_4\nsubseteq G$ and $\delta(G)\ge m$ as desired.  \hfill$\Box$


\end{spacing}
\end{document}